# ROBUST ESTIMATION FOR ARMA MODELS

By Nora Muler,[1] Daniel Peña[2] and Víctor J. Yohai[3]

*Universidad Torcuato di Tella, Universidad Carlos III de Madrid and Universidad de Buenos Aires and CONICET*

This paper introduces a new class of robust estimates for ARMA models. They are M-estimates, but the residuals are computed so the effect of one outlier is limited to the period where it occurs. These estimates are closely related to those based on a robust filter, but they have two important advantages: they are consistent and the asymptotic theory is tractable. We perform a Monte Carlo where we show that these estimates compare favorably with respect to standard M-estimates and to estimates based on a diagnostic procedure.

**1. Introduction.** There are two main approaches to deal with outliers when estimating ARMA models. The first approach is to start estimating the model parameters using maximum likelihood and then analyzing the residuals with a diagnostic procedure to detect outliers. Among others, diagnostic procedures for detecting outliers were proposed by Fox [9], Chang, Tiao and Chen [4], Tsay [23], Peña [22] and Chen and Liu [5]. However diagnostic procedures suffer from the masking problem: when there are several outliers that have similar effects, the outliers may not be detected.

A second approach is to use robust estimates, that is, estimates that are not much influenced by outlying observations. A detailed review of robust procedures for ARMA models can be found in Chapter 8 of Maronna, Martin and Yohai [16]. In that chapter, it is shown that in the case of an AR($p$) model, one outlier at observation $t$ can affect the residuals corresponding to periods $t'$, $t \leq t' \leq t + p$; in the case of an ARMA($p, q$) model with $q > 0$, it can affect all residuals corresponding to periods $t' \geq t$. For this reason estimates based on regular residuals (e.g., M- or S-estimates)

Received June 2006.

[1]Supported in part by a Grant PAV from ANPCYT, Argentina.

[2]Supported in part by MEC Grant SEJ2004-03303, Spain.

[3]Supported in part by Grant X-094 from the Universidad de Buenos Aires, Grant PICT 21407 from ANPCYT, Argentina and Grant PIP 5505 from CONICET, Argentina.

*AMS 2000 subject classifications.* 62F35, 62M10.

*Key words and phrases.* MM-estimates, outliers, time series.

are not very robust. One way to improve the robustness of the estimates is to compute the residuals using the robust filter introduced by Masreliez [20]. This robust filter approximates the one-step predictor in ARMA models with additive outliers. Several authors have proposed estimates that use residuals computed with the Masreliez filter. For instance, Martin, Samarov and Vandaele [19] proposed filtered M-estimates, Martin and Yohai [18] filtered S-estimates and Bianco et al. [1] filtered $\tau$-estimates. However, we can mention two shortcomings of the estimates based on filtered residuals. First, these estimates are asymptotically biased. Second, there is not an asymptotic theory for these estimators, and therefore inference procedures like tests or confidence regions are not available.

In this paper, we propose a new approach to avoid the propagation of the effect of one outlier when computing the innovation residuals of the ARMA model: we define these residuals using an auxiliary model. For this purpose we introduce the *bounded innovation propagation ARMA* (BIP-ARMA) models. With the help of these models, we are able to define estimates for the ARMA model that are highly robust when the series contains outliers.

We show that the mechanisms of the proposed estimates to avoid the propagation of the outliers are similar to those based on robust filters. However, the advantage of these estimates over those based on the robust filters is that they are consistent and asymptotically normal under a perfectly observed ARMA model.

The proposed estimates can be considered as a generalization of the MM-estimates introduced by Yohai [25] for regression. In the first step we define a highly robust residuals scale and in the second step we use a redescending M-estimate, which uses this scale.

For brevity's sake, we have omitted in this paper some of the proofs. All the proofs can be found in Muler, Peña and Yohai [21].

The paper is organized as follows. In Section 2 we introduce the new family of models and show that the corresponding residuals are similar to those obtained with a robust filter. In Section 3 we introduce the proposed estimates. In Section 4 we establish the main asymptotic results: consistency and asymptotic normality. In Section 5 we discuss the computation of the proposed estimates. In Section 6 we discuss robustness properties of the proposed estimates. In Section 7 we present the results of a Monte Carlo study. In Section 8 we show the performance of the different estimates for fitting a monthly real series. In Section 9 we make some concluding remarks. Section 10 is an Appendix with the main proofs of the asymptotic results.

## 2. A new class of bounded nonlinear ARMA models.

2.1. *BIP-ARMA models.* We are going to consider a stationary and invertible ARMA model that can be represented by

$$\phi(B)(x_t - \mu) = \theta(B)a_t, \tag{1}$$

where $a_t$ are i.i.d. random variables with symmetric distribution and where $\phi(B)$ and $\theta(B)$ are polynomial operators given by $\phi(B) = 1 - \sum_{i=1}^{p} \phi_i B^i$ and $\theta(B) = 1 - \sum_{i=1}^{q} \theta_i B^i$ with roots outside the unit circle.

If $a_t$ has first moment we have that $E(x_t) = \mu$. Let $\lambda(B) = \phi^{-1}(B)\theta(B) = 1 + \sum_{i=1}^{\infty} \lambda_i B^i$ and consider the MA($\infty$) representation of the ARMA process

$$x_t = \mu + a_t + \sum_{i=1}^{\infty} \lambda_i a_{t-i}. \tag{2}$$

We can model contaminated ARMA processes with a fraction $\varepsilon$ of outliers by

$$z_t^\varepsilon = (1 - \zeta_t^\varepsilon)x_t + \zeta_t^\varepsilon w_t, \tag{3}$$

where $x_t$ is the ARMA model, $w_t$ is an arbitrary process and $\zeta_t^\varepsilon$ is a process taking values 0 and 1 such that $\lim_{n\to\infty} 1/n(\sum_{i=1}^{n} \zeta_t^\varepsilon) = \varepsilon$. For example, $\zeta_t^\varepsilon$ may be a stationary process such that $E(\zeta_t^\varepsilon) = \varepsilon$. The case of additive outliers corresponds to $w_t = x_t + \nu_t$, where $x_t$ and $v_t$ are independent processes. Replacement outliers correspond to the case that the processes $x_t$ and $w_t$ are independent. According to the dependence structure of the process $\zeta_t^\varepsilon$, we can have additive outliers or patchy outliers. For detail, see Martin and Yohai [17]. Robustness is related to the possibility of accurately estimating the parameter of the central model $x_t$ when we observe the contaminated process $z_t^\varepsilon$.

Another type of outlier are innovation outliers. An ARMA process with innovation outliers occurs when we observe an ARMA process satisfying (1) but the innovations $a_t$ have a heavy-tailed distribution. Regular M-estimates can cope with this type of outlier. See for example Maronna, Martin and Yohai [16].

We will use the following family of auxiliary models

$$y_t = \mu + a_t + \sum_{i=1}^{\infty} \lambda_i \sigma \eta\left(\frac{a_{t-i}}{\sigma}\right), \tag{4}$$

where the $a_t$'s are i.i.d. random variables with symmetric distribution and $\sigma$ is a robust M-scale of $a_t$, which coincides with the standard deviation in the case that the $a_t$'s are normal, the $\lambda_i$'s are the coefficients of $\phi^{-1}(B)\theta(B)$ and $\eta(x)$ is an odd and bounded function. An M-scale of $a_t$ is defined as the solution of the equation $E(\rho(a_t/\sigma)) = b$. We call this model the bounded innovation propagation autoregressive moving average model (BIP-ARMA).

To obtain robust and efficient estimates we will choose $\eta$ bounded and such that there exists $k$ with $\eta(x) = x$ for $|x| \leq k$. More details on how to choose $\rho$, $b$ and $\eta$ are given in Sections 3.1 and 6. Note that in this model the lag effect of a large innovation in period $t$ has a bounded effect on $y_{t+j}$ for any $j \geq 0$ and, since $\lambda_j \to 0$ exponentially when $j \to \infty$, this effect will almost disappear in a few periods.

Note that (4) can also be written as

$$y_t = \mu + a_t - \sigma\eta\left(\frac{a_t}{\sigma}\right) + \sigma\phi^{-1}(B)\theta(B)\eta\left(\frac{a_t}{\sigma}\right)$$

and multiplying both sides by $\phi(B)$, we get

$$\phi(B)y_t = \mu\left(1 - \sum_{i=1}^{p}\phi_i\right) + \phi(B)a_t - \sigma\phi(B)\eta\left(\frac{a_t}{\sigma}\right) + \sigma\theta(B)\eta\left(\frac{a_t}{\sigma}\right)$$

which is equivalent to

$$(5) \quad y_t = a_t + \mu + \sum_{i=1}^{p}\phi_i(y_{t-i} - \mu) - \sum_{i=1}^{r}\left(\phi_i a_{t-i} + (\theta_i - \phi_i)\sigma\eta\left(\frac{a_{t-i}}{\sigma}\right)\right),$$

where $r = \max(p, q)$. If $r > p$, $\phi_{p+1} = \cdots = \phi_r = 0$ and if $r > q$, $\theta_{q+1} = \cdots = \theta_r = 0$.

2.2. *Robust filters and BIP-ARMA models.* Let us analyze the relationship of the BIP-ARMA model and an ARMA model with additive outliers. The BIP-ARMA model can be also be written as $y_t = (1 - \zeta_t^\varepsilon)x_t + \zeta_t^\varepsilon(x_t + \nu_t)$, where $x_t = y_t - a_t + a_t^*$ is an ARMA model satisfying $\phi(B)(x_t - \mu) = \theta(B)a_t^*$, $a_t^* = \sigma\eta(a_t/\sigma)$, $\nu_t = a_t - a_t^*$, $\zeta_t^\varepsilon = I(|a_t| \geq k)$ and $\varepsilon = P(|a_t| \geq k)$. However, in the BIP-ARMA model, $\zeta_t^\varepsilon$ and $\nu_t$ are not independent and they are also not independent of $x_t$.

We will show that the one-step forecast in the BIP-ARMA model is similar to the forecast obtained by using the robust filter for ARMA models introduced by Masreliez [20]. The Masreliez filter was proposed as an approximation to one-step predictors for additive models of the form (3), where $x_t$ is a Gaussian ARMA model, $\zeta_t^\varepsilon$ are i.i.d. Bernoulli variables with $P(\zeta_t^\varepsilon = 1) = \varepsilon$ and $\nu_t$ are i.i.d. normal random variables.

Suppose we have an ARMA series $y_1, \ldots, y_n$ and we suspect that it is contaminated with additive outliers. Assume first that we know the parameters $\phi, \theta, \mu$ and $\sigma$ of the ARMA model. The robust filter computes a "clean" series $y_t^*$, and filtered innovation residuals $\widehat{a}_t$ that are obtained by the following recursive procedure. Suppose the cleaned series $y_1^*, \ldots, y_{t-1}^*$, and the filtered innovation residuals $\widehat{a}_1, \ldots, \widehat{a}_{t-1}$ previous to time $t$ are computed. Since $y_t = \mu - \sum_{i=1}^{\infty} \pi_i(y_{t-i} - \mu) + a_t$, where $\pi(B) = \theta(B)^{-1}\phi(B) = 1 + \sum_{i=1}^{\infty} \pi_i B^i$, the one-step ahead robust forecast of $y_t$ is obtained by replacing the $y_{t-i}$'s by the cleaned values $y_{t-i}^*$'s, that is, the one-step robust forecast of $y_t$ is obtained by

$$\widehat{y}_t^* = \mu - \sum_{i=1}^{\infty} \pi_i(y_{t-i}^* - \mu) = \mu - (\theta(B)^{-1}\phi(B) - 1)(y_t^* - \mu), \tag{6}$$

where $y_t^* = \mu$ for $t \leq 0$. The filtered innovation residual for period $t$ is computed by $\widehat{a}_t^* = y_t - \widehat{y}_t^*$ and the cleaned value $y_t^*$ by

$$y_t^* = \widehat{y}_t^* + s_t \eta^*\left(\frac{\widehat{a}_t^*}{s_t}\right) = y_t - \widehat{a}_t^* + s_t \eta^*\left(\frac{\widehat{a}_t^*}{s_t}\right), \tag{7}$$

where $s_t$ is an estimate of the one-step prediction error scale and where $\eta^*$ has the same properties as those stated for $\eta$, in particular for some $k > 0$ it holds that $\eta^*(u) = u$ for $|u| \leq k$. Observe that if $|\widehat{a}_t^*| \leq k$, then $s_t \eta^*(\widehat{a}_t^*/s_t) = \widehat{a}_t^*$ and $y_t^* = y_t$. Recursive formulae for $s_t$ can be found in Martin, Samarov and Vandaele [19].

We can easily derive from (6) and (7) that

$$\widehat{y}_t^* = \mu + \sum_{i=1}^{p} \phi_i(y_{t-i}^* - \mu) - \sum_{i=1}^{q} \theta_i s_t \eta^*\left(\frac{\widehat{a}_{t-i}^*}{s_t}\right). \tag{8}$$

Now, from (5), the one-step forecast for $y_t$ in the BIP-ARMA model is given by

$$\widehat{y}_t = \mu + \sum_{i=1}^{p} \phi_i\left(y_{t-i} - \mu - a_{t-i} + \sigma\eta\left(\frac{a_{t-i}}{\sigma}\right)\right) - \sum_{i=1}^{q} \theta_i \sigma \eta\left(\frac{a_{t-i}}{\sigma}\right), \tag{9}$$

which is similar to (8) taking as the cleaned series

$$y_t^* = y_t - a_t + \sigma\eta(a_t/\sigma). \tag{10}$$

The main difference is that here $s_t$ is taken constant and equal to $\sigma$. Thus, the filtered residuals used by Martin, Samarov and Vandaele [19] and Bianco et al. [1] are very similar to those of a BIP-ARMA model. In the next section, we will use the model (4) to define robust estimates of the parameters of an ARMA model that may contain additive outliers.

**3. Bounded MM-estimates for ARMA models.** Assume that $y_1, \ldots, y_n$ are observations corresponding to a BIP-ARMA model and that the density of $a_t$ is $f(u)$. The conditional log-likelihood function of $y_{p+1}, \ldots, y_n$ given $y_1, \ldots, y_p$ and the values $a_{p-r+1}^b(\boldsymbol{\beta}, \sigma) = 0, \ldots, a_p^b(\boldsymbol{\beta}, \sigma) = 0$, where $r = \max(p, q)$ can be written as

$$L_c(\boldsymbol{\beta}, \sigma) = \sum_{t=p+1}^{n} \log f(a_t^b(\boldsymbol{\beta}, \sigma)), \tag{11}$$

where from (5), the functions $a_t^b(\boldsymbol{\beta},\sigma)$ are defined recursively for $t \geq p+1$ by

$$
\begin{aligned}
(12) \quad a_t^b(\boldsymbol{\beta},\sigma) &= y_t - \mu - \sum_{i=1}^{p} \phi_i(y_{t-i} - \mu) \\
&\quad + \sum_{i=1}^{r} \left( \phi_i a_{t-i}^b(\boldsymbol{\beta},\sigma) + (\theta_i - \phi_i)\sigma\eta\left(\frac{a_{t-i}^b(\boldsymbol{\beta},\sigma)}{\sigma}\right) \right).
\end{aligned}
$$

In the case of a pure ARMA model, $\eta(u) = u$, (12) reduces to

$$
(13) \quad a_t(\boldsymbol{\beta}) = y_t - \mu - \sum_{i=1}^{p} \phi_i(y_{t-i} - \mu) + \sum_{i=1}^{q} \theta_i a_{t-i}(\boldsymbol{\beta}).
$$

Since ML-estimates are not robust, we will consider M-estimates, which minimize

$$
(14) \quad M_n^b(\boldsymbol{\beta}) = \frac{1}{n-p} \sum_{t=p+1}^{n} \rho\left(\frac{a_t^b(\boldsymbol{\beta},\widehat{\sigma})}{\widehat{\sigma}}\right),
$$

where $\rho$ is a bounded function and $\widehat{\sigma}$ is an estimate of $\sigma$.

We observe that the M-estimates defined in (14) require an estimate $\widehat{\sigma}$ of $\sigma$. This leads us to define in Section 3.2 a two-step procedure for estimating $\boldsymbol{\beta}$ that we call MM-estimates.

3.1. *M-estimates of scale.* Huber [12] introduced the M-estimates of scale. Given a sample $\mathbf{u} = (u_1, \ldots, u_n)$, $u_i \in R$, an M-estimate of scale $S_n(\mathbf{u})$ is defined by any value $s \in (0, \infty)$ satisfying

$$
(15) \quad \frac{1}{n} \sum_{i=1}^{n} \rho\left(\frac{u_i}{s}\right) = b,
$$

where $\rho$ is a function satisfying the following property P1:

P1: $\rho(0) = 0$, $\rho(x) = \rho(-x)$, $\rho(x)$ is continuous, nonconstant and nondecreasing in $|x|$.

We can define two asymptotic breakdown points of an M-estimate of scale: the minimum fraction of outliers that are required to make the estimate infinity, $\epsilon_\infty^*$, and the minimum fraction of inliers that may take the estimate to zero, $\epsilon_0^*$. Huber [13] shows that $\epsilon_\infty^* = b/a$ and $\epsilon_0^* = 1 - b/a$, where $a = \max \rho$. Then, the global breakdown point of these estimates is $\epsilon^* = \min(\epsilon_\infty^*, 1 - \epsilon_\infty^*)$ and taking $b = a/2$, we get a maximum breakdown point of 0.5. To make the M-scale estimate consistent for the standard deviation when the data are normal, we require that $E_\Phi(\rho(x)) = b$ where $\Phi$ is the standard normal distribution.

3.2. *MM-estimates.* The MM-estimates for regression were introduced by Yohai [25] to combine high breakdown points with high efficiency under normal errors. The key idea of the MM-estimates is to compute in the first step a highly robust estimate of the error scale, and in the second step this scale estimate is used to compute an M-estimate of the regression parameters. For time series models, these two steps are not enough to guarantee robustness. This is due to the fact that an outlier in one period not only affects the residual corresponding to this period but it may also affect all the subsequent residuals. To avoid this propagation we define MM-estimates for the ARMA model, where the residuals are computed as in the BIP-ARMA model instead as in the regular ARMA model. Then, the procedure for computing MM-estimates is as follows.

*Step* 1. In this step we obtain an estimate of $\sigma$. For this purpose we consider two estimates of $\sigma$, one using an ARMA model and another using a BIP-ARMA model, and choose the smallest one.

Let $\rho_1$ be a bounded function satisfying P1 and such that if $b = E(\rho_1(u))$, then $b/\max \rho_1(u) = 0.5$. This guarantees that for a normal random sample the M-scale estimator $s$ based on $\rho_1$ converges to the standard deviation and that the breakdown point of $s$ is 0.5. Put

$$\tag{16} \mathcal{B}_{0,\zeta} = \{(\boldsymbol{\phi}, \boldsymbol{\theta}) \in R^{p+q} : |z| \geq 1 + \zeta$$

holds for all the roots $z$ of $\phi(B)$ and $\theta(B)\}$.

Let us call $\mathcal{B}_0 = \mathcal{B}_{0,\zeta}$ for some small $\zeta > 0$ and $\mathcal{B} = \mathcal{B}_{0,\zeta} \times R$. Then, we define an estimate of $\boldsymbol{\beta}$

$$\tag{17} \widehat{\boldsymbol{\beta}}_S = \arg\min_{\boldsymbol{\beta} \in \mathcal{B}} S_n(\mathbf{a}_n(\boldsymbol{\beta}))$$

and the corresponding estimate of $\sigma$

$$\tag{18} s_n = S_n(\mathbf{a}_n(\widehat{\boldsymbol{\beta}}_S)),$$

where $\mathbf{a}_n(\boldsymbol{\beta}) = (a_{p+1}(\boldsymbol{\beta})), \ldots, a_n(\boldsymbol{\beta})$ and $S_n$ is the M-estimate of scale based on $\rho_1$ and $b$.

Let us describe now the estimate corresponding to the BIP-ARMA model. Define $\widehat{\boldsymbol{\beta}}_S^b = (\widehat{\boldsymbol{\phi}}_S^b, \widehat{\boldsymbol{\theta}}_S^b, \widehat{\mu}_s^b)$ by the minimization of $S_n(\mathbf{a}_n^b(\boldsymbol{\beta}, \widehat{\sigma}(\boldsymbol{\phi}, \boldsymbol{\theta})))$ over $\mathcal{B}$. The value $\widehat{\sigma}(\boldsymbol{\phi}, \boldsymbol{\theta})$ is an estimate of $\sigma$ computed as if $(\boldsymbol{\phi}, \boldsymbol{\theta})$ were the true parameters and the $a_t$'s were normal. Then, since in this case the M-scale $\sigma$ coincides with the standard deviation of $a_t$, from (4) we have

$$\sigma^2 = \frac{\sigma_y^2}{1 + \kappa^2 \sum_{i=1}^{\infty} \lambda_i^2(\boldsymbol{\phi}, \boldsymbol{\theta})},$$

where $\kappa^2 = \mathrm{var}(\eta(a_t/\sigma))$ and $\sigma_y^2 = \mathrm{var}(y_t)$. Let $\widehat{\sigma}_y^2$ be a robust estimate of $\sigma_y^2$ and $\kappa^2 = \mathrm{var}(\eta(z))$ where $z$ has $N(0,1)$ distribution. Then, we define

$$\widehat{\sigma}^2(\boldsymbol{\phi}, \boldsymbol{\theta}) = \frac{\widehat{\sigma}_y^2}{1 + \kappa^2 \sum_{i=1}^{\infty} \lambda_i^2(\boldsymbol{\phi}, \boldsymbol{\theta})}. \tag{19}$$

The scale estimate $s_n^b$ corresponding to the BIP-ARMA model is defined by

$$\widehat{\boldsymbol{\beta}}_S^b = (\widehat{\boldsymbol{\phi}}_S^b, \widehat{\boldsymbol{\theta}}_S^b, \widehat{\mu}_s^b) = \arg\min_{\boldsymbol{\beta} \in \mathcal{B}} S_n(\mathbf{a}_n^b(\boldsymbol{\beta}, \widehat{\sigma}(\boldsymbol{\phi}, \boldsymbol{\theta}))) \tag{20}$$

and

$$s_n^b = S_n(\mathbf{a}_n^b(\widehat{\boldsymbol{\beta}}_S^b, \widehat{\sigma}(\widehat{\boldsymbol{\phi}}_S^b, \widehat{\boldsymbol{\theta}}_S^b))), \tag{21}$$

where $\mathbf{a}_n^b(\boldsymbol{\beta}, \sigma) = (a_{p+1}^b(\boldsymbol{\beta}, \sigma), \ldots, a_n^b(\boldsymbol{\beta}, \sigma))$. Our estimate of $\sigma$ is

$$s_n^* = \min(s_n, s_n^b). \tag{22}$$

As we will see in the next section, if the sample is taken from an ARMA model without outliers, asymptotically we obtain $s_n < s_n^b$. We should point out that despite the fact that $\kappa$ was computed as if the $a_t$'s were normal, the asymptotic properties of the estimators are not going to depend on this assumption.

*Step* 2. Consider a bounded function $\rho_2$ such that satisfies P1 and $\rho_2 \leq \rho_1$. This function is chosen so that the corresponding M-estimate is highly efficient under normal innovations. Let

$$M_n(\boldsymbol{\beta}) = \frac{1}{n-p} \sum_{t=p+1}^{n} \rho_2\left(\frac{a_t(\boldsymbol{\beta})}{s_n^*}\right) \tag{23}$$

and

$$M_n^b(\boldsymbol{\beta}) = \frac{1}{n-p} \sum_{t=p+1}^{n} \rho_2\left(\frac{a_t^b(\boldsymbol{\beta}, s_n^*)}{s_n^*}\right). \tag{24}$$

We define the estimates $\widehat{\boldsymbol{\beta}}_M$ and $\widehat{\boldsymbol{\beta}}_M^b$ by the minimization over $\mathcal{B}$ of $M_n(\boldsymbol{\beta})$ and $M_n^b(\boldsymbol{\beta})$, respectively. Then, the MM-estimate $\widehat{\boldsymbol{\beta}}_M^*$ is equal to $\widehat{\boldsymbol{\beta}}_M$ if $M_n(\widehat{\boldsymbol{\beta}}_M) \leq M_n^b(\widehat{\boldsymbol{\beta}}_M^b)$ and is equal to $\widehat{\boldsymbol{\beta}}_M^b$ if $M_n(\widehat{\boldsymbol{\beta}}_M) > M_n^b(\widehat{\boldsymbol{\beta}}_M^b)$.

For instance, we can take $\rho_2(x) = \rho_1(\lambda x)$ with $0 < \lambda < 1$. If $\rho_2''(0) > 0$, $\rho_2$ will be close to a quadratic function when $\lambda$ tends to 0.

REMARK 1. One important problem that will be only briefly mentioned here is that of the robust model selection. One possibility to explore is to adapt to ARMA models the robust finite prediction error (RFPE) selection criterion given in Section 5.12 of Maronna, Martin and Yohai [16] for regression.

In the next section we will show that when the sample is taken from an ARMA model without outliers, for large $n$ the estimate will choose $\widehat{\boldsymbol{\beta}}_M^* = \widehat{\boldsymbol{\beta}}_M$. In our Monte Carlo study of Section 7 we observe that if the sample has enough additive outliers we may have $\widehat{\boldsymbol{\beta}}_M^* = \widehat{\boldsymbol{\beta}}_M^b$. This implies that $\widehat{\boldsymbol{\beta}}_M^*$ and $\widehat{\boldsymbol{\beta}}_M^b$ have the same asymptotic distribution for any $\eta$. However, the efficiency of $\widehat{\boldsymbol{\beta}}_M^b$ for finite sample size depends on $\eta$. If the interval where $\eta$ coincides with the identity increases, the efficiency for finite sample size of $\widehat{\boldsymbol{\beta}}_M^b$ will increase as well, but the propagation of the outliers effect will gain importance and so the estimate will lose robustness.

**4. Asymptotic results.** The main results of this section, stated in Theorems 4 and 6, are the consistency and asymptotic normality of the BMM-estimators for ARMA models. These theorems require to prove first the consistency of S- and the consistency and asymptotic normality of MM-estimators. We stated these results in Theorems 1, 3 and 5, respectively. The link that relates the properties of S- and MM- to those of BMM-estimates are Theorems 2 and 4.

Consider the following assumptions:

P2. The process $y_t$ is a stationary and invertible ARMA$(p,q)$ process with parameter $\boldsymbol{\beta}_0 = (\boldsymbol{\phi}_0, \boldsymbol{\theta}_0, \mu_0) \in \mathcal{B}$ and $E(\log^+|a_t|) < \infty$, where $\log^+ a = \max(\log a, 0)$. The polynomials $\phi_0(B)$ and $\theta_0(B)$ do not have common roots.

P3. The innovation $a_t$ has an absolutely continuous distribution with a symmetric and strictly unimodal density.

P4. $P(a_t \in C) < 1$ for any compact $C$.

P5. The function $\eta$ is continuous, even and bounded.

The following theorem establishes the consistency of the S-estimates based on ARMA models.

THEOREM 1. *Assume that $y_t$ satisfies* P2 *with innovations $a_t$ satisfying* P3. *Assume also that $\rho_1$ is bounded and satisfies* P1 *with $\sup \rho_1 > b$, and that $\psi_1 = \rho_1'$ is bounded and continuous. Then,* (i) *the estimate $\widehat{\boldsymbol{\beta}}_S$ defined in (17) is strongly consistent for $\boldsymbol{\beta}_0$.* (ii) *Let $s_n$ be the scale estimate defined in (18). Then $s_n \longrightarrow s_0$ a.s. where $s_0$ is defined by $E(\rho_1(a_t/s_0)) = b$.*

The next theorem establishes that under a regular ARMA model $\widehat{\boldsymbol{\beta}}_S$ and $\widehat{\boldsymbol{\beta}}_S^b$ are asymptotically equivalent.

THEOREM 2. *Assume that $y_t$ satisfies condition* P2, *with innovations $a_t$ satisfying* P3 *and* P4. *Assume also that $\rho_1$ is bounded and satisfies P1 with $\sup \rho_1 > b$, that $\psi_1 = \rho_1'$ is bounded, continuous and that $\eta$ satisfies* P5. *Then if $y_t$ is not white noise, with probability 1 there exists $n_0$ such that $\widehat{\boldsymbol{\beta}}_S^b = \widehat{\boldsymbol{\beta}}_S$ for all $n \geq n_0$ and then $s_n^*$ defined in (22) verifies $s_n^* \longrightarrow s_0$ a.s.*

Note that when $y_t$ is white noise both models: the regular ARMA and the BIP-ARMA, coincide. Then, the assumption that $y_t$ is not white noise is essential for the validity of the theorem above.

The following theorems shows the consistency of the MM-estimate.

THEOREM 3. *Assume that $y_t$ satisfies condition* P2, *with innovations $a_t$ satisfying* P3. *Assume also that $\rho_i$, $i=1,2$, are bounded and satisfy* P1, $\psi_i = \rho_i'$, $i=1,2$ *are bounded and continuous and that* $\sup \rho_1 > b$. *Then* $\widehat{\boldsymbol{\beta}}_M \longrightarrow \boldsymbol{\beta}_0$ *a.s.*

The next theorem shows that asymptotically under a regular ARMA model $\widehat{\boldsymbol{\beta}}_M$ and $\widehat{\boldsymbol{\beta}}_M^b$ are equivalent.

THEOREM 4. *Suppose that the assumptions of Theorem 3, P4 and P5 hold. Then if $y_t$ is not white noise, with probability 1 there exists $n_0$ such that $\widehat{\boldsymbol{\beta}}_M^b = \widehat{\boldsymbol{\beta}}_M$ for all $n \geq n_0$ and then $\widehat{\boldsymbol{\beta}}_M^* \to \boldsymbol{\beta}_0$ a.s.*

The following theorem shows the asymptotic normality of the MM-estimates.

THEOREM 5. *Suppose that the assumptions of Theorem 3 hold. Moreover, assume that $\psi_2'$ and $\psi_2''$ are continuous and bounded functions and $\sigma_a^2 = E(a_t^2) < \infty$. Then we have*

$$(n-p)^{1/2}(\widehat{\boldsymbol{\beta}}^M - \boldsymbol{\beta_0}) \to_D N(\mathbf{0}, D),$$

*where*

(25)
$$D = \frac{s_0^2 E_{F_0}(\psi_2^2(a_t/s_0))}{E_{F_0}^2(\psi_2'(a_t/s_0))} \begin{pmatrix} \sigma_a^{-2} \mathbf{C}^{-1} & \mathbf{0} \\ \mathbf{0} & \zeta_0^{-2} \end{pmatrix},$$

$$\zeta_0 = -\frac{1 - \sum_{i=1}^p \phi_{0i}}{1 - \sum_{i=1}^q \theta_{0i}}$$

*and $C = (c_{ij})$ is the $(p+q+1) \times (p+q+1)$ matrix given by*

$$c_{i,j} = \sum_{k=0}^{\infty} v_k v_{k+j-i} \qquad \text{if } i \leq j \leq p,$$

$$c_{p+i,p+j} = \sum_{k=0}^{\infty} \varpi_k \varpi_{k+j-i} \qquad \text{if } i \leq j \leq q,$$

$$c_{i,p+j} = -\sum_{k=0}^{\infty} \varpi_k v_{k+j-i} \qquad \text{if } i \leq p, j \leq q, i \leq j,$$

$$c_{i,p+j} = -\sum_{k=0}^{\infty} v_k \varpi_{k+i-j} \qquad \text{if } i \leq p, j \leq q, j \leq i,$$

where $\phi_0^{-1}(B) = 1 + \sum_{i=1}^{\infty} v_i B^i$ and $\theta_0^{-1}(B) = 1 + \sum_{i=1}^{\infty} \varpi_i B^i$. Observe that when the $a_t$'s are normal, $\sigma^2 = \sigma_a^2$.

REMARK 2. When $\rho_2(u) = u^2$, $\widehat{\boldsymbol{\beta}}_M$ is the conditional maximum likelihood estimate corresponding to normal errors. Let $F_0$ be the distribution of $a_t$, then, in this case $s_0^2 E_{F_0}(\psi_2^2(a_t/s_0))/E_{F_0}{}^2(\psi_2'(a_t/s_0)) = \sigma_a^2$. Therefore, the asymptotic efficiency of the MM-estimate with respect to the normal conditional maximum likelihood estimate when the innovations have distribution $F_0$ is

$$\text{EFF}(\psi_2, F_0) = \frac{\sigma_a^2 E_{F_0}^2(\psi_2'(a_t/s_0))}{s_0^2 E_{F_0}(\psi_2^2(a_t/s_0))}. \tag{26}$$

Choosing $\psi_2$ conveniently we can make this efficiency as close to one as desired for the case that $F_0$ is normal.

REMARK 3. The relative efficiency of the MM- and BMM- estimates given by (26) is the same as the one of the M-estimates of location with respect to the mean. This implies the well-known fact that M-estimates are robust for innovation outliers, that is, when $y_t, 1 \leq t \leq n$, correspond to a perfectly observed ARMA model but the distribution $F_0$ of $a_t$ is heavy tailed.

REMARK 4. When $E(a_t^2) = \infty$, the rate of convergence of M-estimates of $\boldsymbol{\phi}$ and $\boldsymbol{\theta}$ may be larger than $n^{-1/2}$, and the asymptotic distribution nonnormal. This case was studied by several authors. See, for example, Davis, Knight and Liu [7] and Davis [8].

REMARK 5. Theorems 1–5 use P3 only to guarantee that for all $\sigma$, the function $g(\mu) = E(\rho((a_t - \mu)/\sigma)$ has a unique minimum at $\mu = 0$. If $g(\mu)$ has a unique minimum at $\overline{\mu} \neq 0$, then the estimates of $\boldsymbol{\phi}$ and $\boldsymbol{\theta}$ are still consistent, but the estimate of $\mu$ will converge to $\mu_0 + \overline{\mu}$.

Finally, from Theorems 4 and 5 we derive the following theorem:

THEOREM 6. *Suppose the assumptions of Theorem 5, P4 and P5 hold. Then $(n-p)^{1/2}(\widehat{\boldsymbol{\beta}}_M^* - \beta_0)$ converges in distribution to a $N(0, D)$ distribution, where $D$ is defined in (25).*

Note that the assumptions of Theorems 2 and 4 include condition P4. However, this condition is not strictly necessary and is included only to simplify the proofs.

All the asymptotic theorems of this section assume that the process is an ARMA model. We conjecture that similar results, consistency and asymptotic normality hold when the observations follows a BIP-ARMA model. The main difficulty to prove these results is to show that the distribution of $a_t^b(\boldsymbol{\beta}, \sigma)$ is asymptotically stationary.

**5. Computation.** We will discuss here how to compute the MM-estimate. We start computing the estimates of Step 1, $\widehat{\boldsymbol{\beta}}_S$ and $\widehat{\boldsymbol{\beta}}_S^b$. According to (15), we can write $S_n^2(\mathbf{a}_n(\boldsymbol{\beta})) = \sum_{t=p+1}^n r_t^2(\boldsymbol{\beta})$, where

$$r_t(\boldsymbol{\beta}) = \frac{S_n(\mathbf{a}_n(\boldsymbol{\beta}))}{(n-p)^{1/2}b^{1/2}} \rho_1^{1/2}\left(\frac{a_t(\boldsymbol{\beta})}{S_n(\mathbf{a}_n(\boldsymbol{\beta}))}\right).$$

Then to compute $\widehat{\boldsymbol{\beta}}_S$ we can use any nonlinear least squares algorithm, for example, a Marquard algorithm. Similarly we can transform the minimization of $S_n(\mathbf{a}_n^b(\boldsymbol{\beta}, \widehat{\sigma}(\boldsymbol{\phi}, \boldsymbol{\theta})))$ in a nonlinear least squares problem. Note that nonlinear least squares algorithms require a good starting point. Since the functions we are minimizing are nonconvex and they may have several local minima, the choice of the starting point is crucial.

If the model has few parameters (e.g., $p+q \le 3$), one way to obtain the starting point is to generate a grid of values of the parameter and choose as initial estimate the one minimizing the objective function. Note that the case of $p+q \le 3$ is very frequent in the case of ARMA applications, where the use of parsimonious models is recommended. Bianco et al. [1] gave an algorithm to compute a highly robust starting point when there are more parameters.

In the second step, to compute $\widehat{\boldsymbol{\beta}}_M$ and $\widehat{\boldsymbol{\beta}}_M^b$ we can use Marquard algorithm using a similar idea and taking as initial estimate the best estimate computed in Step 1.

In our simulations the estimates were defined taking

$$\rho_2(x) = \begin{cases} 0.5x^2, & \text{if } |x| \le 2, \\ 0.002x^8 - 0.052x^6 + 0.432x^4 - 0.972x^2 + 1.792, & \text{if } 2 < |x| \le 3, \\ 3.25, & \text{if } |x| > 3, \end{cases}$$

$\rho_1(x) = \rho_2(x/0.405)$ and $\eta = \rho_2'$. Note that $\rho_1$ and $\rho_2$ are smooth functions that are quadratic in the intervals $(-0.81, 0.81)$ and $(-2, 2)$, respectively. The function $\rho_1$ was chosen so that if we take $b = \max \rho_1/2$ then the scale is consistent to the standard deviation for normal samples. Note that $\eta$ is a redescending function.

For fitting an ARMA$(1,1)$ model to 1000 observations using a MATLAB program, with an initial solution computed with a grid of 20 values in each parameter, the computing time of a BMM-estimate with the choices of $\rho_i$, $i = 1, 2$, and $\eta$ given above is approximately 10 seconds in a personal computer with an AMD Athlon 1.8 GHz processor. For fitting an AR(3) model under the same conditions, the computing time is 1 minute 20 seconds.

**6. Robustness properties.** Several robustness measures can be used for estimates of time-series parameters. Hampel [11] introduced the influence curve to measure the robustness of an estimate under an infinitesimal outlier contamination in the framework of i.i.d. observations. Künsch [14], Martin and Yohai [17] and Mancini, Ronchetti and Trojani [15] give generalizations of the influence curve for estimating time-series parameters. However, because of its infinitesimal character, the influence curve may not be a good measure of the robustness when there is a positive fraction of outlier contamination. For example, it can be proved that for a very small amount of contamination the MM- and BMM-estimates asymptotically coincide and therefore their influence curves also coincide. However, we will see below in this section and in Section 7 that the BMM-estimate is more robust than the MM-estimate. Influence functions for the M-estimates of ARMA models can be found in Martin and Yohai [17].

A more reliable measure of the robustness of an estimate to cope with a positive fraction $\varepsilon$ of contamination is the asymptotic maximum bias. Consider a family of $\varepsilon$-contaminated process

$$z_t^{\varepsilon k} = (1 - \zeta_t^{\varepsilon})x_t + \zeta_t^{\varepsilon} w_t^k \tag{27}$$

as in (3) where $k \in K$ and $(x_t, \zeta_t^{\varepsilon}, w_t^k)$ is stationary. Suppose also the distribution of the uncontaminated process $x_t$ depends on a parameter $\boldsymbol{\gamma} \in \Gamma \subset R^j$. As an example, we can consider the family of additive outliers models, which is obtained taking $w_t^k = x_t + k$, with $k \in R$.

Suppose that for each sample size $n$ we have an estimate $\widehat{\boldsymbol{\gamma}}_n$ of $\boldsymbol{\gamma}$ and let $\widehat{\boldsymbol{\gamma}}_\infty(L) = (\widehat{\gamma}_{\infty 1}(L), \ldots, \widehat{\gamma}_{\infty j}(L))$ be the almost sure limit of $\widehat{\boldsymbol{\gamma}}_n$ when applied to a process with distribution $L$. The bias of the $i$th component $\widehat{\boldsymbol{\gamma}}_\infty$ when

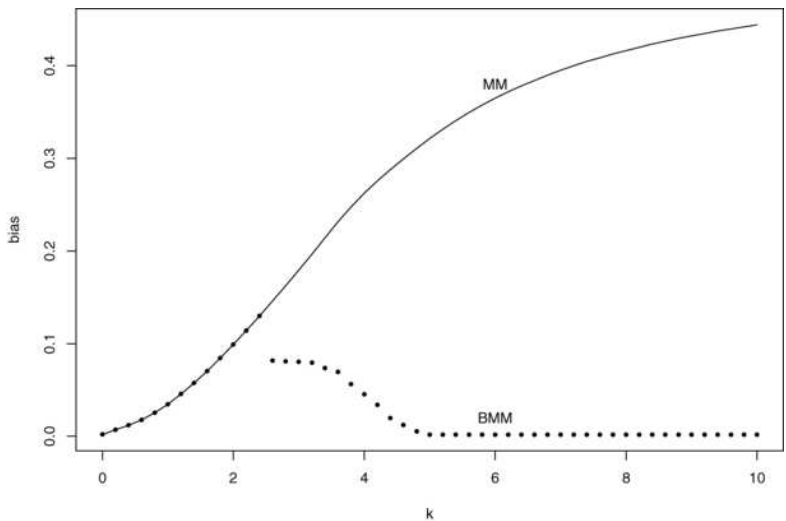

FIG. 1. *Bias curve for the* AR(1) *model with* $\phi = 0.5$ *and 10% of additive outliers where* $k$ *is the outlier size.*

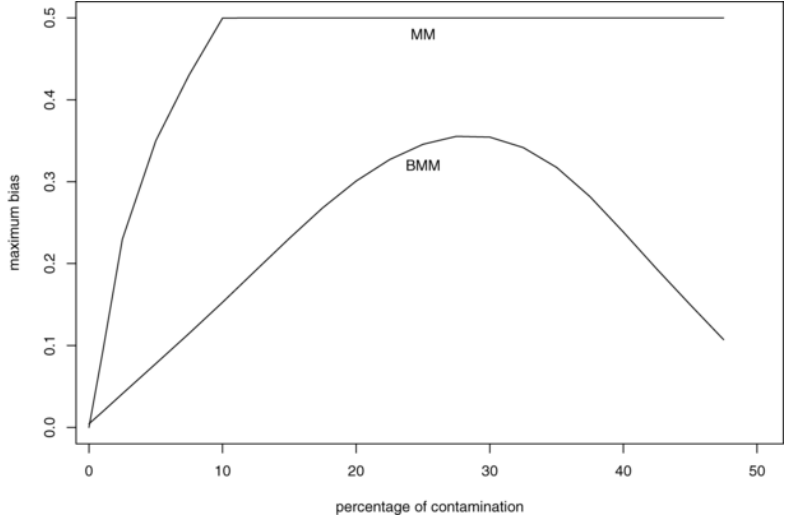

Fig. 2. *Maximum bias for the* AR(1) *model with* $\phi = 0.5$.

applied to $z_t^{\varepsilon k}$ as defined in (27) is

$$B(\widehat{\gamma}_{\infty i}, \gamma, \varepsilon, k) = |\widehat{\gamma}_{\infty i}(\mathcal{L}(z_t^{\varepsilon k})) - \gamma_i|,$$

where $\mathcal{L}(z_t^{\varepsilon k})$ denotes the distribution of the process $z_t^{\varepsilon k}$. The maximum asymptotic bias of the $i$th component is defined by

$$MB(\widehat{\gamma}_{\infty i}, \gamma, \varepsilon) = \sup_{k \in K} B(\widehat{\gamma}_{\infty i}, \gamma, \varepsilon, k).$$

We have approximately computed the maximum bias curves of the MM- and BMM-estimates for Gaussian AR(1) and MA(1) models with additive outliers ($w_t^k = x_t + k$) and where the $\zeta_t^\varepsilon$ are i.i.d. Bernoulli variables. To simplify the computation we eliminate the intercept from these models by assuming it to be known and null. The asymptotic value of the estimate is approximated using samples of size 10000. We found that for samples size larger than 10000 the changes in the estimate are negligible.

In Figure 1, we show the bias curves of the MM- and BMM-estimates for the AR(1) model with $\phi = 0.5$ and $\varepsilon = 0.1$. In Figure 2, we show the maximum biases curves for the MM- and BMM-estimates under the same model. In Figure 3, we show the maximum bias curve for the BMM-estimate under a MA(1) model with parameter $\theta = -0.5$.

In both cases, we observe that the BMM-estimate has a smaller maximum bias than the MM-estimate. We also observe that the behavior of the MM- is different from the BMM-estimate. After the contamination is larger than some value $\varepsilon^*$ the maximum bias of the MM-estimate is constantly equal to the value of the estimated parameter. This means that the asymptotic value of the estimate becomes 0 independently of the true value of the parameter. This value $\varepsilon^*$ corresponds to the breakdown point notion proposed by Genton and Lucas [10]. For the AR(1) model the value $\varepsilon^*$ depends on $\phi$. For the MA(1) model $\varepsilon^* = 0$. Instead, the behavior of the BMM-estimate is different and apparently the estimate does not break down. A surprising feature of its maximum bias curve is that for very large $\varepsilon$ the maximum bias starts decreasing. This can be explained as follows: when $\varepsilon$ is large, the probability of obtaining a patch of two or more outliers increases. The effect of a patch of outliers is to increase the correlation of the series. Therefore, in the case of the AR(1) model with $\phi$ positive and MA(1) with $\theta$ negative it prevents that the parameter further approximates to zero for outliers with fixed size $k$. We also computed the maximum bias curves for other values of parameters $\phi$ and $\theta$ and the results were similar.

We conjecture that the robust behavior of the BMM-estimate under additive outlier contamination also holds when we observe any contaminated process $z_t^\varepsilon$ as given in (3). The reason is that since the BIP-ARMA model includes a built-in filtering to compute the residuals, a small fraction of outliers will affect only a small fraction of residuals. Therefore, since the loss function of the BMM-estimate is bounded, the estimate will not be largely affected by a small fraction of large residuals. We compute maximum bias curves for the case of replacement outliers ($w_t^k = k$), obtaining similar results as for the case of additive outliers.

**7. A Monte Carlo study.** We have performed a Monte Carlo study to compare several estimates for ARMA models. We have simulated three

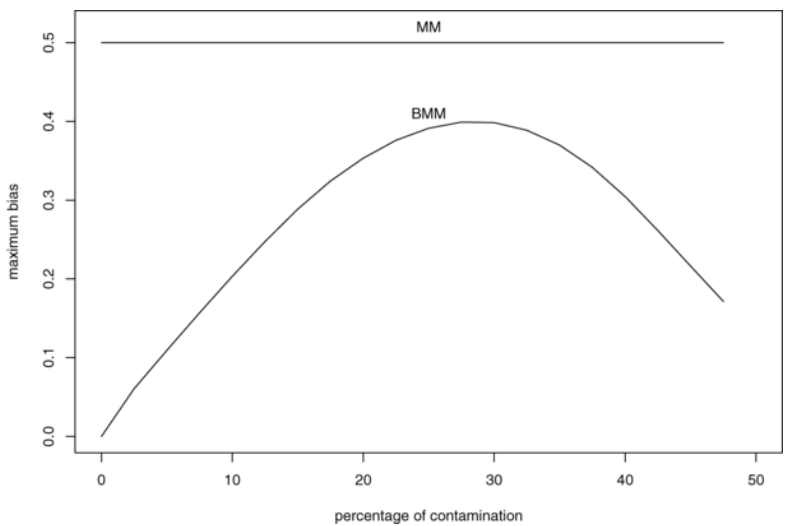

Fig. 3. *Maximum bias of the BMM for the* MA(1) *model with* $\theta = 0.5$.

Gaussian stationary ARMA models considering the case that the series do not contain outliers and the case that the series have 10% of equally spaced in time additive outliers. The size of the additive outliers was taken equal to 4 and 6. The sample size in the simulations is 200 and the Monte Carlo study was done with 500 replications.

The estimates considered in this study are (i) the normal conditional maximum likelihood estimate (MLE), (ii) an MM-estimate where the residuals are computed as in a regular ARMA model (MM), (iii) an MM-estimate where the residuals are compared with the ones of a BIP-ARMA model (BMM), (iv) an estimate based on the diagnostic procedure described in Chen and Liu [5]. The cutoff point for outlier rejection is chosen by the Bonferroni inequality as $c = \Phi^{-1}(1 - (0.05/n))$, where $\Phi$ is the $N(0,1)$ distribution function. We denote this estimate by ($\text{CTC}_B$). (v) The same as in (iv) but the cutoff point is $c = 3$ ($\text{CTC}_3$). (vi) The tau filtered estimate proposed by Bianco et al. [1]. We denote this estimate by (FTAU).

The estimates MM and BMM are based on the functions $\rho_1$ and $\rho_2$ and $\eta$ described in Section 5. We have simulated three models: AR(1) with $\phi = 0.5$, MA(1) with $\theta = -0.5$ and ARMA(1,1) with $\phi = 0.5$ and $\theta = -0.5$. The value of $\mu$ was zero for all these models. However, since the estimates are shift equivariant, their behavior does not depend on the value of this parameter. For brevity's sake, we only show here (in Tables 1 and 2) the results for the AR(1) and the MA(1) model. The results for the ARMA(1,1) models are similar and can be found in Muler, Peña and Yohai [21].

The relative efficiency with respect to the MLE when there are no outliers varies in the case of $\phi$ and $\theta$ from 80% to 86% for the estimate BMM, from 80% to 91% for the estimate MM, is almost 100% for the CTC estimates and varies from 65% to 67% for the FTAU. The efficiency of all the estimates of $\mu$ is very high. We also observe that under additive outlier contamination, the estimate BMM of $\phi$ and $\theta$ behaves much better than those corresponding to the estimates MM, $\text{CTC}_B$ and $\text{CTC}_3$. The performance of the estimates FTAU and BMM are comparable.

The errors of the MSEs shown on Tables 1 and 2 are smaller than 15% with probability 0.95. However, since all the estimates were computed with the same samples, the errors of the differences between the MSEs of any two estimates are much smaller, thus making comparisons possible.

**8. An example.** This example deals with a monthly series of inward movement of residential telephone extensions in a fixed geographic area from January 1966 to May 1973 (RESEX). The series was analyzed by Brubacher [2] and by Martin, Samarov and Vandaele [19], who identified an AR(2) model for the differenced series $y_t = x_t - x_{t-12}$, where $x_t$ is the observed series.

Table 1

*MSE for the AR(1) model with $\phi = 0.5$ without outliers and with 10% of equally spaced additive outliers*

|  | No outliers | | Out. size 4 | | Out. size 6 | |
| --- | --- | --- | --- | --- | --- | --- |
| **Estimate** | $\mu$ | $\phi = 0.5$ | $\mu$ | $\phi = 0.5$ | $\mu$ | $\phi = 0.5$ |
| MLE | 0.017 | 0.0036 | 0.189 | 0.103 | 0.394 | 0.189 |
| MM | 0.018 | 0.0045 | 0.024 | 0.085 | 0.028 | 0.132 |
| BMM | 0.018 | 0.0042 | 0.021 | 0.014 | 0.019 | 0.0048 |
| $CTC_B$ | 0.017 | 0.0036 | 0.185 | 0.103 | 0.364 | 0.189 |
| $CTC_3$ | 0.017 | 0.0036 | 0.148 | 0.096 | 0.057 | 0.047 |
| FTAU | 0.019 | 0.0054 | 0.032 | 0.011 | 0.028 | 0.0076 |

Table 3 displays the value of the estimates MLE, MM, BMM, $CTC_3$ and the FTAU together with the MAD-scale of the residuals. We can see that the estimated values of the parameters of the MLE and the $CTC_3$ are quite different from the robust estimates MM, BMM and FTAU. The estimate $CTC_B$ gives the same result as $CTC_3$ (it detects the same outliers) and is omitted from the table.

Table 2

*MSE for the MA(1) model with $\theta = -0.5$ without outliers and with 10% of equally spaced additive outliers*

|  | No outliers | | Out. size 4 | | Out. size 6 | |
| --- | --- | --- | --- | --- | --- | --- |
| **Estimate** | $\mu$ | $\theta$ | $\mu$ | $\theta$ | $\mu$ | $\theta$ |
| MLE | 0.010 | 0.0042 | 0.178 | 0.128 | 0.380 | 0.215 |
| MM | 0.012 | 0.0046 | 0.015 | 0.115 | 0.016 | 0.159 |
| BMM | 0.012 | 0.0052 | 0.015 | 0.025 | 0.012 | 0.0065 |
| $CTC_B$ | 0.011 | 0.0042 | 0.174 | 0.130 | 0.345 | 0.218 |
| $CTC_3$ | 0.011 | 0.0044 | 0.136 | 0.125 | 0.042 | 0.064 |
| FTAU | 0.012 | 0.0065 | 0.020 | 0.031 | 0.017 | 0.021 |

Table 3

*Estimates of the parameters of the RESX series*

| **Estimates** | $\mu$ | $\phi_1$ | $\phi_2$ | **MAD** |
| --- | --- | --- | --- | --- |
| MLE | 2.69 | 0.48 | $-0.17$ | 1.70 |
| MM | 1.18 | 0.34 | 0.31 | 1.43 |
| BMM | 1.74 | 0.42 | 0.36 | 1.24 |
| $CTC_3$ | 3.44 | 1.14 | $-0.74$ | 1.86 |
| FTAU | 1.71 | 0.27 | 0.49 | 1.10 |

Figure 4 shows the data $y_t$ obtained differentiating the observed data as $y_t = x_t - x_{t-12}$ and the cleaned values as in (10), which are seen to be almost coincident except at outlier locations.

**9. Concluding remarks.** We have presented two families of estimates for ARMA models: MM-estimates and BMM-estimates. The BMM-estimates use a mechanism that avoids the propagation of the full effect of the outliers to the subsequent residual innovations. To make this mechanism compatible with consistency when the true model is ARMA, we consider two estimates: one is obtained fitting a regular ARMA model and the other fitting a BIP-ARMA model, where the propagation of the effect of outliers is bounded. Then, the estimate that fits better to the data is selected. We have shown in Sections 6 and 7 that, at least for additive outliers, the BMM-estimates are much more robust than the MM-estimates and quite comparable with the FTAU-estimates. The main advantage of the BMM-estimates over the FTAU-estimates is that an asymptotic theory is now available and this makes inference with BMM-estimates possible. The Monte Carlo results of Section 7 also show that the BMM-estimate compares favorably with the estimate based on the Chen and Liu [5] diagnostic procedure.

## APPENDIX

Suppose that we have the infinite sequence of observations $\mathbf{Y}_t = (\ldots, y_{t-k}, \ldots, y_{t-1}, y_t)$ generated by a stationary and invertible $\text{ARMA}(p,q)$ process up to time $t$ with parameter $\boldsymbol{\beta}_0$. Given any $\boldsymbol{\beta} = (\boldsymbol{\phi}, \boldsymbol{\theta}, \mu)$ such that the polynomials $\phi(B)$ and $\theta(B)$ have all the roots outside the unit circle, let us define $a_t^e(\boldsymbol{\beta}) = \theta^{-1}(B)\phi(B)(y_t - \mu)$. Then $a_t^e(\boldsymbol{\beta}_0) = a_t$ and $a_t^e(\boldsymbol{\beta})$'s satisfy the following recursive relationship

$$a_t^e(\boldsymbol{\beta}) = y_t - \mu - \sum_{i=1}^{p} \phi_i(y_{t-i} - \mu) + \sum_{i=1}^{q} \theta_i a_{t-i}^e(\boldsymbol{\beta}).$$

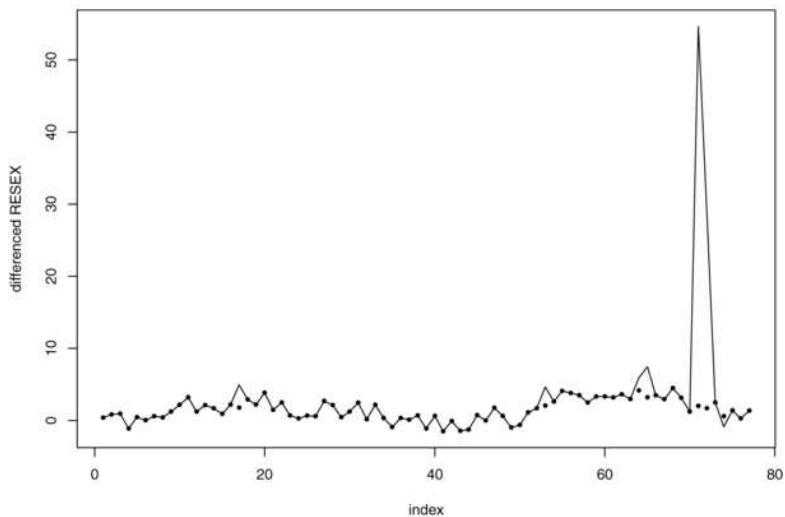

FIG. 4. *Differenced RESEX Series: observed (solid line) and filtered (dots) values.*

In the case that $a_t$ has finite first moment, we have that $a_t^e(\boldsymbol{\beta}) = y_t - E(y_t|\mathbf{Y}_{t-1})$, where the conditional expectation is taken assuming that the true value of the parameter vector is $\boldsymbol{\beta}$.

We will use the following notation. Given a function $g(\mathbf{u}): R^k \to R$, we define $\nabla g(\mathbf{u})$ as the column vector of dimension $k$ whose $i$th element is $\nabla_i g(\mathbf{u}) = \partial g(\mathbf{u})/\partial u_i$ and $\nabla^2 g(\mathbf{u})$ is the $k \times k$ matrix whose $(i,j)$ element is $\nabla_{ij}^2 g(\mathbf{u}) = \partial^2 g(\mathbf{u})/\partial u_i \, \partial u_j$.

The next lemma proved in Lemma 2 of Muler, Peña and Yohai [21] gives the Fisher Consistency of the S-estimate when we have all the past observations.

LEMMA 1. *Assume that $y_t$ satisfies condition* P2 *with innovations satisfying* P3. *Assume that $\rho_1$ is a bounded function satisfying condition* P1, *define $s(\boldsymbol{\beta})$ by $E(\rho_1(a_t^e(\boldsymbol{\beta})/s(\boldsymbol{\beta}))) = b$. Then, if $\boldsymbol{\beta} \in \mathcal{B}$ and $\boldsymbol{\beta} \neq \boldsymbol{\beta}_0$ we have $s_0 = s(\boldsymbol{\beta}_0) < s(\boldsymbol{\beta})$.*

The proofs of the next two lemmas can be found in Lemmas 5 and 6 of Muler, Peña and Yohai [21].

LEMMA 2. *Under the assumptions of Theorem 1, for any $d > 0$ we have,*

$$\lim_{n \to \infty} \sup_{\boldsymbol{\beta} \in \mathcal{B}_0 \times [-d,d]} |S_n(\mathbf{a}_n(\boldsymbol{\beta})) - s(\boldsymbol{\beta})| = 0 \quad a.s.$$

LEMMA 3. *Under the assumptions of Theorem 1, there exists $d > 0$ satisfying*

$$\liminf_{n \to \infty} \inf_{|\mu| > d, (\boldsymbol{\phi}, \boldsymbol{\theta}) \in \mathcal{B}_0} S_n(\mathbf{a}_n(\boldsymbol{\beta})) > s_0 + 1 \quad a.s.$$

PROOF OF THEOREM 1. Take $\varepsilon > 0$ arbitrarily small and let $d$ be as in Lemma 3. By the dominated convergence theorem it is easy to show that $s(\boldsymbol{\beta})$ is continuous. Then by Lemma 1, there exists $0 < \gamma < 1$ such that

$$\min_{\boldsymbol{\beta} \in \mathcal{B}_0 \times [-d,d], \|\boldsymbol{\beta} - \boldsymbol{\beta}_9\| \geq \varepsilon} s(\boldsymbol{\beta}) \geq s_0 + \gamma.$$

Therefore by Lemma 2, there exist $n_1$ such that for $n \geq n_1$

$$\min_{\boldsymbol{\beta} \in \mathcal{B}_0 \times [-d,d], \|\boldsymbol{\beta} - \boldsymbol{\beta}_0\| \geq \varepsilon} S_n(\boldsymbol{\beta}) \geq s_0 + \gamma/2$$

and $S_n(\boldsymbol{\beta}_0) \leq s_0 + \gamma/4$. Moreover by Lemma 3, there exists $n_2$ such that for $n \geq n_2$

$$\inf_{|\mu|>d,(\boldsymbol{\phi},\boldsymbol{\theta})\in\mathcal{B}_0} S_n(\mathbf{a}_n(\boldsymbol{\beta})) > s_0 + \gamma \qquad \text{a.s.}$$

Therefore, for $n \geq \max(n_1, n_2)$ it holds that $\|\widehat{\boldsymbol{\beta}}_S - \boldsymbol{\beta}_0\| < \varepsilon$ and this proves the theorem. $\square$

The next three lemmas will be used to prove Theorem 2.

The proof of the next two lemmas can be found in Lemmas 7 and 8 of Muler, Peña and Yohai [21].

LEMMA 4. *Assume that $y_t$ satisfies condition* P2. *Given $d > 0$ and $\widetilde{\sigma} > 0$, there exist constants $C > 0$ and $0 < \nu < 1$ such that*

$$\sup_{\boldsymbol{\beta}\in\mathcal{B}_0\times[-d,d]} \sup_{0<\sigma\leq\widetilde{\sigma}} |a_t^b(\boldsymbol{\beta},\sigma) - y_t| \leq C\left(\widetilde{\sigma} + \nu^t \sum_{i=1}^p |y_i|\right), \qquad t \geq p+1.$$

LEMMA 5. *Under the assumptions of Theorem 2, given $d > 0$, there exists $\delta > 0$ such that*

$$\liminf_{n\to\infty} \inf_{\boldsymbol{\beta}\in\mathcal{B}_0\times[-d,d]} S_n(\mathbf{a}_n^b(\boldsymbol{\beta}, \widehat{\boldsymbol{\sigma}}(\boldsymbol{\phi},\boldsymbol{\theta}))) > s_0 + \delta \qquad a.s.$$

The proof of the next lemma can be found in Lemma 9 of Muler, Peña and Yohai [21].

LEMMA 6. *Under the assumptions of Theorem 2, there exists $d > 0$ such that,*

$$\lim\inf_{n\to\infty} \inf_{|\mu|>d,(\boldsymbol{\phi},\boldsymbol{\theta})\in\mathcal{B}_0} S_n(\mathbf{a}_n^b(\boldsymbol{\beta}, \widehat{\sigma}(\boldsymbol{\phi},\boldsymbol{\theta}))) > s_0 + 1 \qquad a.s.$$

PROOF OF THEOREM 2. From Lemmas 5 and 6 we have that there exists $\delta > 0$ such that

$$\lim\inf_{n\to\infty} \inf_{\boldsymbol{\beta}\in\mathcal{B}} S_n(\mathbf{a}_n^b(\boldsymbol{\beta}, \widehat{\sigma}(\boldsymbol{\phi},\boldsymbol{\theta}))) > s_0 + \delta \qquad \text{a.s.}$$

But, by Theorem 1(ii) we have that $\widehat{\boldsymbol{\beta}}_S$ satisfies $\lim_{n\to\infty} S_n(\mathbf{a}_n(\widehat{\boldsymbol{\beta}}_S)) = s_0$ a.s. This proves the theorem. $\square$

The following four lemmas will be used to prove Theorem 3. The proofs can be found in Lemmas 10, 11, 12 and 13 of Muler, Peña and Yohai [21].

LEMMA 7. *Assume that $y_t$ satisfies condition* P2 *with innovations satisfying* P3 *and assume that $\rho_2$ satisfies condition* P1 *with $\rho_2$ bounded. Let us call $m(\boldsymbol{\beta}) = E(\rho_2(a_t^e(\boldsymbol{\beta})/s_0))$, then $\boldsymbol{\beta}_0 = \arg\min_{\boldsymbol{\beta}\in\mathcal{B}} m(\boldsymbol{\beta})$.*

LEMMA 8. *Assume that $y_t$ satisfies condition P2 and $\rho_2$ condition P1. Define*

$$M_n^e(\boldsymbol{\beta}) = \frac{1}{n-p} \sum_{t=p+1}^{n} \rho_2\left(\frac{a_t^e(\boldsymbol{\beta})}{s_n^*}\right). \tag{28}$$

*Then, we have,*

$$\lim_{n\to\infty} \sup_{\boldsymbol{\beta}\in\mathcal{B}_0\times[-d,d]} \left| M_n^e(\boldsymbol{\beta}) - E\left(\rho_2\left(\frac{a_t^e(\boldsymbol{\beta})}{s_0}\right)\right)\right| = 0 \quad \text{a.s.}$$

*for all $d > 0$.*

LEMMA 9. *Under the assumptions of Theorem 3, we have*

$$\lim_{n\to\infty} \sup_{\boldsymbol{\beta}\in\mathcal{B}_0\times[-d,d]} |M_n(\boldsymbol{\beta}) - M_n^e(\boldsymbol{\beta})| = 0 \quad \text{a.s.}$$

LEMMA 10. *Under the assumptions of Theorem 3, there exists $d > 0$ and $\delta > 0$ such that*

$$\liminf_{n\to\infty} \inf_{|\mu|>d,(\boldsymbol{\phi},\boldsymbol{\theta})\in\mathcal{B}_0} M_n(\boldsymbol{\beta}) \geq m(\boldsymbol{\beta}_0) + \delta \quad \text{a.s.,}$$

*where $m(\boldsymbol{\beta}_0)$ is defined in Lemma 7.*

PROOF OF THEOREM 3. Follows from Lemmas 8–10 using similar arguments as those used in the proof of Theorem 1. □

The next two lemmas will be used to prove Theorem 4.

LEMMA 11. *Under the assumptions of Theorem 3, for all $d > 0$ there exists $\delta > 0$ such that*

$$\liminf_{n\to\infty} \inf_{\boldsymbol{\beta}\in\mathcal{B}_0\times[-d,d]} M_n^b(\boldsymbol{\beta}) \geq m(\boldsymbol{\beta}_0) + \delta \quad \text{a.s.,}$$

*where $m(\boldsymbol{\beta}_0)$ is defined in Lemma 7.*

PROOF. It is similar to the proof of Lemma 5. □

The proof of the next lemma can be found in Lemma 15 of Muler, Peña and Yohai [21].

LEMMA 12. *Under the assumptions of Theorem 3, there exists $d > 0$ and $\delta > 0$ such that*

$$\liminf_{n\to\infty} \inf_{|\mu|>d,(\boldsymbol{\phi},\boldsymbol{\theta})\in\mathcal{B}_0} M_n^b(\boldsymbol{\beta}) \geq m(\boldsymbol{\beta}_0) + \delta \quad \text{a.s.,}$$

*where $m(\boldsymbol{\beta}_0)$ is defined as in Lemma 7.*

PROOF OF THEOREM 4. From Lemmas 11 and 12 we have that there exists $\delta > 0$ such that $\liminf_{n\to\infty} \inf_{\boldsymbol{\beta}\in\mathcal{B}} M_n^b(\boldsymbol{\beta}) \geq m(\boldsymbol{\beta}_0) + \delta$. Theorem 3 implies that $\lim_{n\to\infty} M_n(\widehat{\boldsymbol{\beta}}_M) = m(\boldsymbol{\beta}_0)$ a.s. This proves the theorem. $\square$

The next five lemmas will be used to prove Theorem 5.

LEMMA 13. *Under the assumptions of Theorem 5, we have*

$$\frac{1}{(n-p)^{1/2}} \sum_{t=p+1}^{n} \nabla\rho_2\left(\frac{a_t^e(\boldsymbol{\beta}_0)}{s_0}\right) \to_D N(\mathbf{0}, V_0),$$

*where*

(29) $$V_0 = E\left(\nabla\rho_2\left(\frac{a_t^e(\boldsymbol{\beta}_0)}{s_0}\right)\nabla\rho_2\left(\frac{a_t^e(\boldsymbol{\beta}_0)}{s_0}\right)'\right).$$

PROOF. This lemma follows immediately from the Central Limit Theorem for Martingales (see Theorem 24.3, Davidson [6]). For details, see Lemma 16 of Muler, Peña and Yohai [21]. $\square$

LEMMA 14. *Under the assumptions of Theorem 5 we have*

$$\lim_{n\to\infty} \frac{1}{(n-p)^{1/2}} \left\| \sum_{t=p+1}^{n} \left(\nabla\rho_2\left(\frac{a_t^e(\boldsymbol{\beta}_0)}{s_n^*}\right) - \nabla\rho_2\left(\frac{a_t^e(\boldsymbol{\beta}_0)}{s_0}\right)\right) \right\| = 0$$

*in probability.*

PROOF. The proof is similar to the one of Lemma 5.1 in Yohai [24] for MM-estimates in the case of regression. The details can be seen in Lemma 17 of Muler, Peña and Yohai [21]. $\square$

The proof of the next lemma can be found in Lemma 18 of Muler, Peña and Yohai [21].

LEMMA 15. *Under the assumptions of Theorem 5, we have*

$$\lim_{n\to\infty} \frac{1}{(n-p)^{1/2}} \left\| \sum_{t=p+1}^{n} \left(\nabla\rho_2\left(\frac{a_t(\boldsymbol{\beta}_0)}{s_n^*}\right) - \nabla\rho_2\left(\frac{a_t^e(\boldsymbol{\beta}_0)}{s_n^*}\right)\right) \right\| = 0 \qquad a.s.$$

The proof of the next lemma can be found in Lemma 19 of Muler, Peña and Yohai [21].

LEMMA 16. *Under the assumptions of Theorem 5 we have for all $d > 0$,*

(i)
$$\lim_{n\to\infty} \sup_{\boldsymbol{\beta}\in\mathcal{B}_0\times[-d,d]} \left\| \frac{1}{n-p} \sum_{t=p+1}^{n} \nabla^2 \rho_2\left(\frac{a_t^e(\boldsymbol{\beta})}{s_n^*}\right) - E\left(\nabla^2 \rho_2\left(\frac{a_t^e(\boldsymbol{\beta})}{s_0}\right)\right) \right\| = 0$$
*a.s.,*

*where $\|A\|$ denotes the $l_2$ norm of the matrix $A$.*

(ii)
$$E\left(\nabla^2 \rho_2\left(\frac{a_t^e(\boldsymbol{\beta}_0)}{s_0}\right)\right) = \frac{1}{s_0^2} E\left(\psi_2'\left(\frac{a_t}{s_0}\right)\right) E(\nabla a_t^e(\boldsymbol{\beta}_0) \nabla a_t^e(\boldsymbol{\beta}_0)')$$

*and this matrix is nonsingular.*

The proof of the next lemma can be found in Lemma 20 of Muler, Peña and Yohai [21].

LEMMA 17. *Under the assumptions of Theorem 5, we have,*

$$\lim_{n\to\infty} \sup_{\boldsymbol{\beta}\in\mathcal{B}_0\times[-d,d]} \frac{1}{n-p} \left\| \sum_{t=p+1}^{n} \left(\nabla^2 \rho_2\left(\frac{a_t(\boldsymbol{\beta})}{s_n^*}\right) - \nabla^2 \rho_2\left(\frac{a_t^e(\boldsymbol{\beta})}{s_n^*}\right)\right) \right\| = 0 \quad a.s.$$

*for all $d > 0$.*

PROOF OF THEOREM 5. The estimate $\widehat{\boldsymbol{\beta}}_M$ satisfies $\sum_{t=p+1}^{n} \nabla \rho_2(a_t(\widehat{\boldsymbol{\beta}}_M)/s_n^*) = 0$. Then, using the Mean Value Theorem we have

$$(30) \quad \sum_{t=p+1}^{n} \nabla \rho_2\left(\frac{a_t(\boldsymbol{\beta}_0)}{s_n^*}\right) + \left(\sum_{t=p+1}^{n} \nabla^2 \rho_2\left(\frac{a_t(\boldsymbol{\beta}^*)}{s_n^*}\right)\right)(\widehat{\boldsymbol{\beta}}_M - \boldsymbol{\beta}_0) = 0,$$

where $\boldsymbol{\beta}^*$ is an intermediate point between $\widehat{\boldsymbol{\beta}}_M$ and $\boldsymbol{\beta}_0$.

From Theorem 3 we have that $\widehat{\boldsymbol{\beta}}_M \to \boldsymbol{\beta}_0$ a.s. Take $d > 0$ so that $d > |\mu_0|$, then with probability one there exists $n_0$ such that $\widehat{\boldsymbol{\beta}}_M \in \mathcal{B}_0 \times [-d,d]$ for all $n \geq n_0$. From Lemmas 16(i) and 17 we get

$$(31) \quad \lim_{n\to\infty} \sup_{\boldsymbol{\beta}\in\mathcal{B}_0\times[-d,d]} \left\| \frac{1}{n-p} \sum_{t=p+1}^{n} \left(\nabla^2 \rho_2\left(\frac{a_t(\boldsymbol{\beta})}{s_n^*}\right) - E\left(\nabla^2 \rho_2\left(\frac{a_t^e(\boldsymbol{\beta})}{s_0}\right)\right)\right) \right\| = 0 \quad \text{a.s.}$$

Put

$$(32) \quad A_n = \frac{1}{n-p} \sum_{t=p+1}^{n} \nabla^2 \rho_2\left(\frac{a_t(\boldsymbol{\beta}^*)}{s_n^*}\right).$$

Then, since $\boldsymbol{\beta}^* \to \boldsymbol{\beta}_0$ a.s. and $E(\nabla^2 \rho_2(a_t^e(\boldsymbol{\beta})/s_0))$ is continuous in $\boldsymbol{\beta}$ we have that

$$(33) \qquad \lim_{n \to \infty} A_n = E\left(\nabla^2 \rho_2\left(\frac{a_t^e(\boldsymbol{\beta}_0)}{s_0}\right)\right) \qquad \text{a.s.}$$

Therefore from Lemma 16(ii), for $n$ large enough $A_n$ is nonsingular. Then, from (30) for large enough $n$ we have $(n-p)^{1/2}(\widehat{\boldsymbol{\beta}}_M - \boldsymbol{\beta}_0) = A_n^{-1} \mathbf{c}_n$, where

$$\mathbf{c}_n = \frac{1}{(n-p)^{1/2}} \sum_{t=p+1}^{n} \nabla \rho_2\left(\frac{a_t(\boldsymbol{\beta}_0)}{s_n^*}\right).$$

From Lemmas 13, 14 and 15 we have that $\mathbf{c}_n \to_D N(\mathbf{0}, V_0)$ and then from (33), and we get $(n-p)^{1/2}(\widehat{\boldsymbol{\beta}}_M - \boldsymbol{\beta}_0) \to_D N(\mathbf{0}, V_1^{-1} V_0 V_1^{-1})$, where $V_1 = E(\nabla^2 \rho_2(a_t^e(\boldsymbol{\beta}_0)/s_0))$.

In Theorem 5 of Muler, Peña and Yohai [21] it is proved that

$$(34) \qquad \nabla \rho_2\left(\frac{a_t^e(\boldsymbol{\beta}_0)}{s_0}\right) = \frac{\psi_2(a_t/s_0)}{s_0} \mathbf{v}_t,$$

where $\mathbf{v}_t$ is the stationary process vector of dimension $(p+q+1)$ defined by

$$\mathbf{v}_{tj} = \begin{cases} -\phi_0^{-1}(B) a_{t-j}, & \text{if } 1 \le j \le p, \\ \theta_0^{-1}(B) a_{t-j-p}, & \text{if } p+1 \le j \le p+q, \\ \zeta_0, & \text{if } j = p+q+1, \end{cases}$$

where $\zeta_0 = -(1 - \sum_{j=1}^{p} \phi_{0j})/(1 - \sum_{j=1}^{p} \theta_{0j})$. Then

$$(35) \qquad E\left(\nabla \rho_2\left(\frac{a_t^e(\boldsymbol{\beta}_0)}{s_0}\right) \nabla \rho_2\left(\frac{a_t^e(\boldsymbol{\beta}_0)}{s_0}\right)'\right) = \frac{E(\psi_2^2(a_t/s_0))}{s_0^2} E(\mathbf{v}_t \mathbf{v}_t').$$

Differentiating $\nabla \rho(a_t^e(\boldsymbol{\beta})/s_0)$ we obtain

$$(36) \qquad \nabla^2 \rho_2\left(\frac{a_t^e(\boldsymbol{\beta}_0)}{s_0}\right) = \frac{1}{s_0^2} \psi_2'\left(\frac{a_t}{s_0}\right) \mathbf{v}_t \mathbf{v}_t' + \frac{1}{s_0} \psi_2\left(\frac{a_t}{s_0}\right) \nabla^2 a_t^e(\boldsymbol{\beta}_0).$$

Since $\nabla^2 a_t^e(\boldsymbol{\beta}_0)$ is independent of $a_t$ and $E(\psi_2(a_t/s_0)) = 0$ we have $E(\psi_2(a_t/s_0) \nabla^2 a_t^e(\boldsymbol{\beta}_0)) = 0$, and then from (36), since $a_t$ and $v_t$ are independent we get

$$(37) \qquad E\left(\nabla^2 \rho_2\left(\frac{a_t(\boldsymbol{\beta}_0)}{s_0}\right)\right) = \frac{1}{s_0^2} E\left(\psi_2'\left(\frac{a_t}{s_0}\right)\right) E(\mathbf{v}_t \mathbf{v}_t').$$

Hence, from (35) and (37) we obtain

$$V_1^{-1} V_0 V_1^{-1} = s_0^2 \frac{E(\psi_2^2(a_t/s_0))}{E(\psi_2'(a_t/s_0))^2} E(\mathbf{v}_t \mathbf{v}_t')^{-1}.$$

Finally, it is straightforward (see, e.g., Bustos and Yohai [3]) to show that

$$E(\mathbf{v}_t \mathbf{v}_t') = \begin{pmatrix} \sigma_a^2 \mathbf{C} & \mathbf{0} \\ \mathbf{0} & \zeta_0^2 \end{pmatrix},$$

where $C$ is defined in the statement of Theorem 5 and $\sigma_a^2 = E(a_t^2)$. □

N. MULER
DEPARTAMENTO DE MATHEMATICAS
 Y ESTADÍSTICA
UNIVERSIDAD TORCUATA DI TELLA
MINONES 2159
1428 BUENOS AIRES
ARGENTINA
E-MAIL: nmuler@utdt.edu

D. PEÑA
DEPARTAMENTO DE ESTADÍSTICA
UNIVERSIDAD CARLOS III DE MADRID
MADRID 126, 28903 GETAFE
MADRID
SPAIN
E-MAIL: dp.rector@uc3m.es

V. J. YOHAI
CIUDAD UNIVERSITARIA
PABELLÓN 1
BUENOS AIRES 1428
ARGENTINA
E-MAIL: vyohai@dm.uba.ar